\documentclass[11pt,twoside]{article}
\usepackage{amsmath,amssymb,latexsym}

\usepackage[english]{babel}
\usepackage{amsmath}
\usepackage{pst-all}

\usepackage{theorem}
\newtheorem{theorem}{Theorem}[section]
\newtheorem{remark}[theorem]{Remark}
\newtheorem{definition}[theorem]{Definition}

\title{The largest small Polytopes}
\author{Andreas Klein and Markus Wessler}
\begin{document}
\maketitle

\begin{abstract}
The aim of this paper is the determination of the largest
$n$-dimensional polytope with $n+3$ vertices of unit diameter.
This is a special case of a more general problem Graham proposes in
\cite{Graham:1975}.
\end{abstract}

\section{Introduction}
We know that among the geometric objects in Euclidean $n$-space
with given diameter the sphere has maximal volume.
A natural question arises if we consider polytopes instead, and in fact
this problem has been considered several times.
In this paper we deal with the following question:
Given natural numbers $k$ and $n$,
which polytope with $k$ vertices of unit diameter
in Euclidean $n$-space has the largest volume?
To this end we define the following volume function:

\begin{definition}
Let $n \geq 2$ and $k \geq n+1$
be positive integers. Then we define $V(n,k)$ to be
the maximum volume of a polytope with $k$ vertices of unit diameter
in Euclidean $n$-space.
\end{definition}

Let us briefly recall the following well-known results.
For $k$ odd, Reinhardt showed \cite{Reinhardt:1922}
that $V(2,k)$ is achieved by the plain regular $k$-gon.
In this case we have
$$V(2, k) = \frac{k}{2} \cdot \cos \left(\frac{\pi}{k}\right)
\cdot \tan \left(\frac{\pi}{2k}\right).$$
It is also known that $V(2,4)$ is achieved by the square
(though not in a unique way, see \cite{Schaeffer:1958}).

However, for even $k > 4$, it is definitely not the plain regular $k$-gon
which has maximal area.
In fact, Graham showed \cite{Graham:1975}
that for $k = 6$ the largest area is not
obtained by the regular hexagon. The latter one is
approximately equal to 0.6495, whereas $V(2,6) = 0.6749 \ldots$.

Considering an $n$-simplex, we observe:
\begin{remark}\label{tetra}
For every $n \geq 2$, $V(n , n+1)$ is achieved by the regular $n$-simplex,
and
we have
$$V(n , n+1) = \frac{1}{n!}\sqrt{\frac{n+1}{2^n}}.$$
\end{remark}

We proceed as follows. For $k = n+2$ the optimal configuration is
obtained in a similar way, as we shall see in the next
section.
In order to determine $V(n, n+3)$ we first consider the special case
$V(3,6)$ which in fact gives rise to a more general procedure leading
to the main result of this paper, Theorem \ref{main}.
We end with an outlook over open problems and some concluding remarks.

\section{The calculation of $V(n,n+2)$}

Each configuration of $n+2$ points in $n$-dimensional space is
topologically equivalent to a double pyramid, which we obtain in the
following way. Let us choose $n$ points in a hyperplane $H$ forming an
$(n-1)$-simplex $S$ and two points $P_1$ and $P_2$ lying on opposite
sides of $S$.  Then we have
$$V(\text{double pyramid}) = \frac{1}{n} \cdot \text{height} \cdot
\text{base},$$ where the height is, of course, bounded by $d(P_1 ,
P_2)$, hence by 1, and the base is bounded
by $V(n-1, n)$. On the other hand, this maximum is really achieved
(though not in unique way). Thus we obtain the following
\begin{theorem}\label{receq}
For $n \ge 2$ we have
$$V(n,n+2) = \frac{1}{n} \cdot V(n-1, n)$$
(where $V(1,2)=1$). Hence, by \ref{tetra},
$$V(n,n+2) = \frac{1}{n!} \sqrt{\frac{n}{2^{n-1}}}.$$
\end{theorem}

\section{The calculation of $V(3,6)$}
The construction of the maximal $n$-dimensional polytopes with $n+3$
vertices is best understood if we look at the
special case $n=3$ first.

Without loss of generality, we may restrict our investigation to
polyhedra with triangular faces. By Euler's polyhedra formula we
obtain that each polyhedron with $6$ vertices must have $8$ faces and
$12$ edges. There are two topologically distinct cases:
\begin{enumerate}
\item No vertex has a valence greater than $4$. In this case the
  polyhedron is topologically equivalent to the regular
  octahedron. (We shall call this the octahedral case.)
\item At least one vertex has valence $5$. In this case the polyhedron
  is topologically equivalent to a pyramid over a pentagon which is
partitioned into three triangles.
In particular, there are two vertices with valence
  $5$. (We shall call this the pyramidal case.)
\end{enumerate}
\begin{figure}[ht]
  \begin{center}
    \begin{pspicture}(0,0)(3.5,3)
      \psline(2,0)(2.4,1.2)
      \psline(2,0)(3.4,1.4)
      \psline(2,0)(0.4,1.6)
      \psline[linestyle=dashed](2,0)(1.4,1.8)
      \psline(0.4,1.6)(2.4,1.2)(3.4,1.4)
      \psline[linestyle=dashed](0.4,1.6)(1.4,1.8)(3.4,1.4)
      \psline(2,3)(2.4,1.2)
      \psline(2,3)(3.4,1.4)
      \psline(2,3)(0.4,1.6)
      \psline[linestyle=dashed](2,3)(1.4,1.8)
    \end{pspicture}
    \hspace{1cm}
    \begin{pspicture}(0,0)(3.5,3)
      \psline(2,0)(2.4,1.2)
      \psline(2,0)(3.5,1.4)
      \psline(2,0)(1.6,1.3)
      \psline(2,0)(0.7,1.8)

      \psline(0.7,1.8)(1.6,1.3)(2.4,1.2)(3.5,1.4)

      \psline(2,3)(2.4,1.2)
      \psline(2,3)(3.5,1.4)
      \psline(2,3)(1.6,1.3)
      \psline(2,3)(0.7,1.8)

      \psline[linestyle=dashed](2,0)(2,3)
    \end{pspicture}
    \caption{The octahedral and the pyramidal case}
  \end{center}
\end{figure}
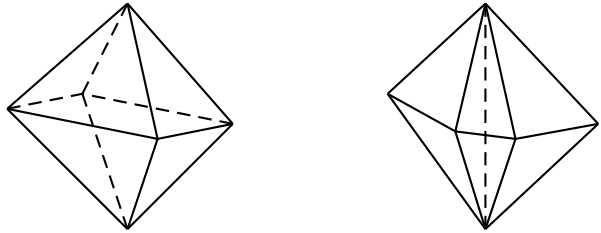

The next step is to determine the maximal volume in each of these
cases. It will turn out that the volume achieved in the octahedral
case is much smaller than the volume achieved in the pyramidal case.

\subsection{The octahedral case} \label{octa}
We fix two opposite vertices $P$ and $Q$. Since each vertex has valence
$4$, we know that for all the other vertices $P_1,P_2,P_3,P_4$ the line
segments $P_iP$ and $P_iQ$ are edges of the polyhedron.

Let $P'_1,P'_2,P'_3,P'_4$ be the projections of $P_1,P_2,P_3,P_4$ to
the plane $p$ that is perpendicular to $PQ$ and that intersects
$PQ$ in its center.
It is obvious that $d(P'_i,P'_j) \le d(P_i,P_j)$ and
$\max(d(P'_i,P),d(P'_i,Q)) \le \max(d(P_i,P),d(P_i,Q))$. Thus the
polyhedron with vertices $Q$, $P$, $P'_1$, $P'_2$, $P'_3$ and $P'_4$
has diameter less than $1$.

Since the volume of $PQP'_iP'_k$ equals the volume of $PQP_iP_k$ (the
height and the area of the base is the same in both pyramids), the
volume
of $QPP'_1P'_2P'_3P'_4$ equals the
volume of $QPP_1P_2P_3P_4$. Thus we can restrict ourselves to the case
that the vertices $P_1$, $P_2$, $P_3$ and $P_4$ lie in the plane
$p$.

But in this case we know that the volume of the polytope is bounded by
$\frac{1}{3}V(2,4) = \frac{1}{6}$. (In case of the regular octahedron
we can achieve equality.)

\subsection{The pyramidal case} \label{pyra}
Now we investigate the pyramidal case. Let $P$ and $Q$ be the vertices
with valence $5$ and $P_1$, $P_2$, $P_3$ and $P_4$ the other vertices
of the polyhedron. Let $p$ be the plane perpendicular to $PQ$
intersecting $PQ$ in its center. As in the octahedral case we conclude
that the volume of the polyhedron does not change if we project the
points $P_1$, $P_2$, $P_3$ and $P_4$ into the plane $p$.
Thus we may again assume that the vertices $P_1$, $P_2$,
$P_3$ and $P_4$ lie in $p$.

Let $P_5 = PQ \cap p$ and $h=d(P,Q)$. Since $d(P,P_i) \le 1$
$(i=1,\dots 4)$ we obtain $d(P_5,P_i)^2+\frac{h^2}{4} \le 1$,
i.e.\ $d(P_5,P_i) \le r = \sqrt{1-\frac{h^2}{4}}$. Since $h \in
[0,1]$, we obtain $r \in [\frac{\sqrt{3}}{2},1]$.

The volume of the polytope is $\frac{1}{3}hS$ where $S$ is the area of
the pentagon $P_1P_2P_3P_4P_5$. Thus we must solve the planar problem
to maximize the area of $P_1P_2P_3P_4P_5$ depending on $r$.

The pentagon $P_1P_2P_3P_4P_5$ satisfies $d(P_5,P_i) \le r$ for
$i=1, \dots, 4$ and $d(P_i,P_j) \le 1$ for $1 \le i < j \le 4$. We
define the \emph{diameter graph} $D$ of the pentagon by:
\begin{itemize}
\item The vertices of $D$ are the points $P_1,\dots,P_5$.
\item $\{P_i,P_5\}$ ($1\le i \le 4$) is an edge of $D$ if and only if
$d(P_i,P_5)=r$.
\item $\{P_i,P_j\}$ ($1\le i < j \le 4$) is an edge of $D$
if and only if $d(P_i,P_j)=1$.
\end{itemize}
In the following we identify an edge $\{P_i,P_j\}$ of $D$ with the
line segment $P_iP_j$.

As in \cite{Graham:1975} (Fact 2) we conclude that the diameter
graph of the pentagon with maximal area is connected.

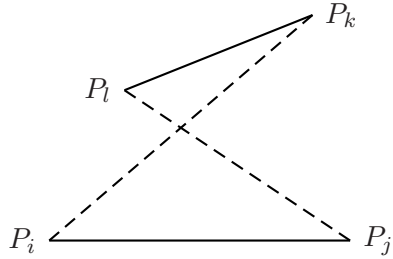
\begin{figure}[ht]
\begin{center}
\begin{pspicture}(0,0)(4,3)
  \psline(0,0)(4,0)
  \psline(1,2)(3.5,3)
  \psline[linestyle=dashed](0,0)(3.5,3)
  \psline[linestyle=dashed](4,0)(1,2)
  \uput[180](0,0){$P_i$}
  \uput[0](4,0){$P_j$}
  \uput[180](1,2){$P_l$}
  \uput[0](3.5,3){$P_k$}
\end{pspicture}
\caption{$D$ has a linear trackleation}\label{triangle}
\end{center}
\end{figure}

Suppose $P_iP_j$ and $P_kP_l$ have no point
in common, then the triangle inequation yields
$d(P_i,P_j)+d(P_k,P_l) < d(P_i,P_k) + d(P_l,P_j)$ (see Figure
\ref{triangle}). Thus $\{P_i, P_j\}$ or $\{P_k, P_l\}$ is not an
edge of $D$.

Graphs with this property are said (by Conway) to have a linear
trackleation. By a result of Woodall \cite{Woodall:1971} the graph $D$
must be one of the following:

\begin{figure}[ht]
\begin{center}
\begin{tabular}{ccc}
(a)
\begin{pspicture}(0,0)(2,2)
  \psline{*-*}(0,1)(1.7,2)
  \psline{*-*}(0,1)(2,1.3)
  \psline{*-*}(0,1)(2,0.7)
  \psline{*-*}(0,1)(1.7,0)
\end{pspicture}
&
(b)
\begin{pspicture}(0,0)(2,2)
  \psline{*-*}(0,1)(1.7,2)
  \psline{*-*}(0,1)(2,1.3)
  \psline{*-*}(0,1)(2,0.7)
  \psline{*-*}(0,1)(1.7,0)
  \psline{*-*}(1.7,2)(1.7,0)
\end{pspicture}
&
(c)
\begin{pspicture}(0,0)(2,2)
  \psline{*-*}(0,1)(1.7,2)
  \psline{*-*}(0,1)(2,1)
  \psline{*-*}(0,1)(1.7,0)
  \psline{*-*}(1.7,2)(1.7,0)
  \psline{*-*}(1.7,2)(0.4,0.1)
\end{pspicture}
\\ \\ \\
(d)
\begin{pspicture}(0,0)(2,2)
  \psline{*-*}(0,2)(2,2)
  \psline{*-*}(0.3,1)(2,2)
  \psline{*-*}(0.5,0.7)(2,2)
  \psline{*-*}(0,2)(1.5,0.7)
\end{pspicture}
&
(e)
\begin{pspicture}(0,0)(2,2)
  \psline{*-*}(1,2)(0.5,0)
  \psline{*-*}(1,2)(1.5,0)
  \psline{*-*}(2,1.5)(0.5,0)
  \psline{*-*}(0,1.5)(1.5,0)
\end{pspicture}
&
(f)
\begin{pspicture}(0,0)(2,2)
  \psline{*-*}(1,2)(0.5,0)
  \psline{*-*}(1,2)(1.5,0)
  \psline{*-*}(2,1.5)(0.5,0)
  \psline{*-*}(0,1.5)(1.5,0)
  \psline{*-*}(2,1.5)(0,1.5)
\end{pspicture}
\end{tabular}
\caption{Possible linear trackleations of $D$}
\end{center}
\end{figure}
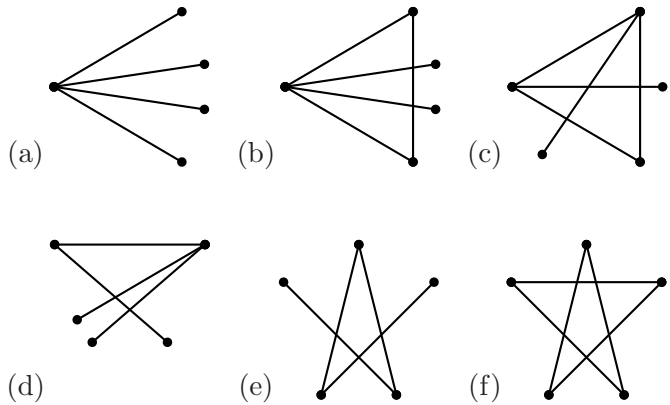

The problem of determining the area of the pentagon is now reduced to
an examination of each of the six cases.

An easy calculation reveals that in the cases (a)-(d) the pentagon
has an area less then $0.567$ for each possible value of $r$ in
$[\frac{\sqrt{3}}{2},1]$.
(For example, in case (b) the pentagon lies in a sixth part of a
circle with radius $1$ and therefore the polygon has an area less than
$0.53$.) As we shall see later, the maximal area in case (f) is
always larger than $0.58$.

Thus the only remaining cases are (e) and (f). We shall prove that in
case (e) there is no \emph{local} maximum and therefore the maximal
area of the pentagon is obtained in case (f).

\begin{figure}[ht]
  \begin{center}
    \begin{pspicture}(0,0)(3,2)
      \psline(0.5,1)(2.1,0)(1.5,2)(0.9,0)(2.7,1.2)
      \pspolygon[linestyle=dashed](0.5,1)(1.5,2)(2.7,1.2)(2.1,0)(0.9,0)
      \uput{0.05}[0](1.1,0.65){\tiny $\alpha_1$}
      \uput{0.15}[180](2,0.7){\tiny $\alpha_2$}
      \uput{0.5}[270](1.5,2){\tiny $\beta$}
    \end{pspicture}
    \caption{The case (e)} \label{case-e}
  \end{center}
\end{figure}
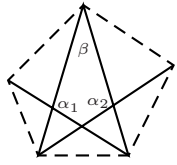

In case (e) we have the situation shown in Figure \ref{case-e}.

We must maximize the area of pentagon $P_1P_2P_3P_4P_5$ depending
on $\alpha_1,\alpha_2$ and $\beta$. By elementary geometry it follows
immediately that a local maximum can only be achieved for
$\alpha_1=\alpha_2=\frac{\pi}{2}$. But
then it is clear that the maximal area is achieved for the case
$\beta=0$. This means that there is no local maximum in case (e).

\begin{figure}[ht]
  \begin{center}
    a)
    \begin{pspicture}(-1,0)(3,2)
      \pspolygon(1.3,0)(0,0.5)(1.7,0.8)(0.7,0)(1,1.6)
      \pspolygon[linestyle=dashed](1.3,0)(0.7,0)(0,0.5)(1,1.6)(1.7,0.8)
      \uput[0](1.7,0.8){$P_1$}
      \uput[0](1.3,0){$P_2$}
      \uput[180](0.7,0){$P_3$}
      \uput[180](0,0.5){$P_4$}
      \uput[90](1,1.6){$P_5$}
    \end{pspicture}
    \hspace{1cm}
    b)
    \begin{pspicture}(-2,0)(4,2)
      \pspolygon(1.3,0)(0,0.6)(2,0.6)(0.7,0)(1,1.6)
      \pspolygon[linestyle=dashed](1.3,0)(0.7,0)(0,0.6)(1,1.6)(2,0.6)
      \uput[0](2,0.6){$P_1=(\frac{1}{2},y)$}
      \uput[0](1.3,0){$P_2=(x,0)$}
      \uput[180](0.7,0){$P_3=(-x,0)$}
      \uput[180](0,0.6){$P_4=(-\frac{1}{2},y)$}
      \uput[90](1,1.6){$P_5=(0,z)$}
    \end{pspicture}
  \end{center}
  \caption{The case (f)} \label{5gon}
\end{figure}
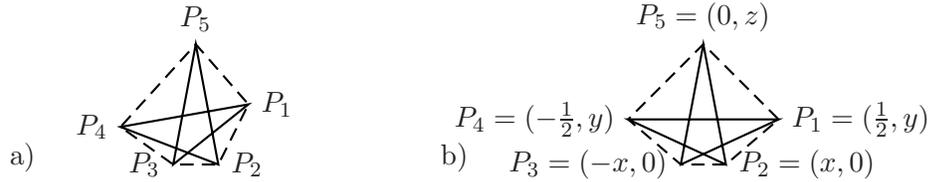

Thus we are left with the last case (f). This is shown in
Figure \ref{5gon}.a.

A simple but tedious
calculation shows that in order to maximize the area of the pentagon,
it is necessary that $\angle{P_4P_2P_3}=\angle{P_2P_3P_1}$. Thus we are
left with the
symmetric case shown in Figure \ref{5gon}.b. We obtain
$y=\sqrt{1-(\frac{1}{2}+x)^2}$ and $z=\sqrt{r^2-x^2}$. The area of the
pentagon is
$$ A(r,x) = \frac{1}{2}x\sqrt{3-4x-4x^{2}}+\frac{1}{2}\sqrt
{{r}^{2}-{x}^{2}}.$$
This expression has a unique local maximum for $x \in [0,\frac{1}{2}]$
which can be found by setting the first derivative equal to zero. The
resulting optimal value $x_0(r)$ is a solution of a sixth order algebraic
equation. Since $r \ge \frac{\sqrt{3}}{2}$ we obtain that
$A(r,x_0(r)) >
A\left(\frac{\sqrt{3}}{2},x_0\left(\frac{\sqrt{3}}{2}\right)\right)
=0.5862\ldots$. This
proves that the maximal area of the pentagon is obtained in case
(f) and not in any of the cases (a)-(e).

Finally, the maximal volume of the polyhedron can be found by maximizing
the
expression
$\frac{1}{3}h
A\left(\sqrt{1-\frac{h^2}{4}},x_0\left(\sqrt{1-\frac{h^2}{4}}\right)\right)$
for
$h \in [0,1]$. Some more calculations reveal that the maximum is
obtained for $h=1$. In this case the volume is $$V(3,6)=0.1954\ldots$$
and there exists a unique polyhedron that archive this maximum.

\section{The calculation of $V(n,n+3)$}
Now we are ready to generalize the arguments of the previous section to
higher dimensions.

As in the $3$-dimensional case we have two topologically distinct
possibilities:

In the octahedral case we have $6$ vertices with valence $n+1$
and $n-3$ vertices of valence $n+2$. Let $P$ and $Q$ be two vertices of
valence $n+1$ and $R_1,\dots,R_{n-3}$ the vertices with valence
$n+2$. Let $p$ be the plane orthogonal to $PQR_1\dots R_{n-3}$
intersecting the $(n-2)$-simplex $PQR_1\dots R_{n-3}$ in the center of
its surrounding sphere. We can generalize the projection argument of
section
\ref{octa} to see that in this case the volume of the polytope is less
than $\frac{1}{n}V(2,4)V(n-2,n-1)$.

In the pyramidal case we find $n-1$ vertices $P_1,\dots,P_{n-1}$ with
valence $n-2$. Let $p$ be the plane orthogonal to
$P_1\dots P_{n-1}$ intersecting the $(n-2)$-simplex
$P_1\dots P_{n-1}$ in the center of the surrounding sphere. If we
project the remaining four points to $p$ we obtain the planar
optimization problem of section \ref{pyra}. In the $n$-dimensional
case,
$r \in \left[\sqrt{1-\frac{(n-2)}{2(n-1)}},1\right]$. (The distance of
the center of
$P_1\dots P_{n-1}$ to the vertices is at most
$\frac{n-2}{n-1}\sqrt{\frac{n-1}{2(n-2)}}$.)
Since now $r$ can be smaller than in section \ref{pyra}, we must
improve the bounds in the cases (a)-(d), but nevertheless, we find
that the area of the pentagon is still maximal in case (f).
Thus we can proceed in the same way as in the $3$-dimensional case and
conclude that the maximal volume is achieved if the $(n-2)$-simplex
$P_1\dots P_{n-1}$ has maximal volume.

\begin{theorem}\label{main}
For $n \ge 3$ let $r=\sqrt{1-\frac{(n-2)}{2(n-1)}}$.
Then we have
$$V(n,n+3) = \frac{A(r,x_0(r))}{n}V(n-2,n-1),$$
where $A$ is the function defined in section \ref{pyra} with local
maximum at $x_0(r)$.
\end{theorem}

In particular, for $n \to \infty$ we have $r \to \frac{1}{\sqrt{2}}$
and $A(r,x_0(r)) \to 0.5002\ldots$, hence

\begin{equation}
  \lim_{n \to \infty} n\frac{V(n,n+3)}{V(n-2,n-1)} = 0.5002\ldots
  \label{limeq2}.
\end{equation}

We remark that
\begin{equation}\label{limeq1}
  \lim_{n \to \infty} n \cdot \frac{V(n, n+2)}{V(n-1, n)} = 1
\end{equation}
(see Theorem \ref{receq}).

\section{Concluding Remarks}
In principle, the preceding techniques (consideration of the
topologically distinct cases, projection, linear trackleation) may be
applied to the determination of $V(n, n+k)$ for $k > 3$. However,
Bender and Wormald \cite{Bender:1988} showed that
$$\frac{1}{972(k-1)(2k-5)(3k-6)}\binom{4k-10}{k+2}$$
is a good approximation for the number of topologically distinct cases
in three dimensions.
The number of possible linear trackleations of the $(2n)$-gon is
$$ \frac{1}{8m} \sum_{\stackrel{d|m}{\text{\tiny$d$ odd}}}
\phi(d)4^{m/d}+4^{m-2}+2^{m-1}-1.$$
Thus the number of topologically distinct cases and the number of
possible linear trackleations in each case, grows exponentially.

However, one could ask if, corresponding to the limit formulae
(\ref{limeq2}) and (\ref{limeq1}), we can determine
$$\lim_{n \to \infty} n\frac{V(n,n+k)}{V(n-k+1,n-k+2)}$$
for $k>3$.

We remark that the maximal polytopes with $n+1$, $n+2$ or
$n+3$ vertices have an axis of symmetry. The problem
whether each maximal polytope has at least one axis of symmetry (see
\cite{Graham:1975}) is still open.

The calculations in section \ref{pyra} suggest the following
generalized problem. Given a symmetric $(k \times k)$-matrix $D$,
determine the largest $k$-gon $P_1\dots P_k$ with
$d(P_i,P_j) \le D_{i,j}$.

\vspace{1cm}

\scriptsize

Andreas Klein\\
Universit\"at Kassel\\
Fachbereich 17 (Mathematik und Informatik)\\
D-34109 Kassel\\
\tt klein@mathematik.uni-kassel.de

Markus Wessler\\
Universit\"at Kassel\\
Fachbereich 17 (Mathematik und Informatik)\\
D-34109 Kassel\\
\tt wessler@mathematik.uni-kassel.de

\end{document}